\documentclass[10pt]{article}
\usepackage{amsmath,amsthm,amsfonts}

\input epsf.tex
\newdimen\epsfxsize
\newdimen\epsfysize
\def\qed{\vrule height5pt width3pt depth.5pt}

\theoremstyle{plain}
\newtheorem{thm}{Theorem}[section]
\newtheorem{cor}[thm]{Corollary}

\newtheorem{conj}{Conjecture}[section]
\newtheorem{rem}{Remark}[section]

\begin{document}

\title{Bounds on Mosaic Knots}

% infomation for author
\author{J. Alan Alewine \\ H. A.  Dye \\
David Etheridge \\ Irina Gardu\~{n}o \\ Amber Ramos  }

\maketitle

\begin{abstract}
We investigate relationships between bounds on the crossing number and the mosaic number of mosaic knots.
\end{abstract}

\section{Introduction}
Lomonaco and Kauffman introduced mosaic knots in \cite{lou-sam}. We first review the defintion of a mosaic knot and then investigate bounds on the crossing number and the mosaic number of a knot.

\begin{figure}[htb] \epsfysize = 0.7 in
\centerline{\epsffile{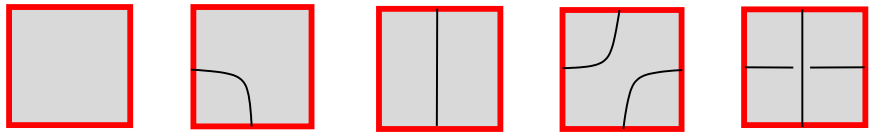}}
\caption{The mosaic tiles, up to rotations of 90 degrees }
\label{fig:tiles}
\end{figure}

\begin{figure}[htb] \epsfysize = 0.7 in
\centerline{\epsffile{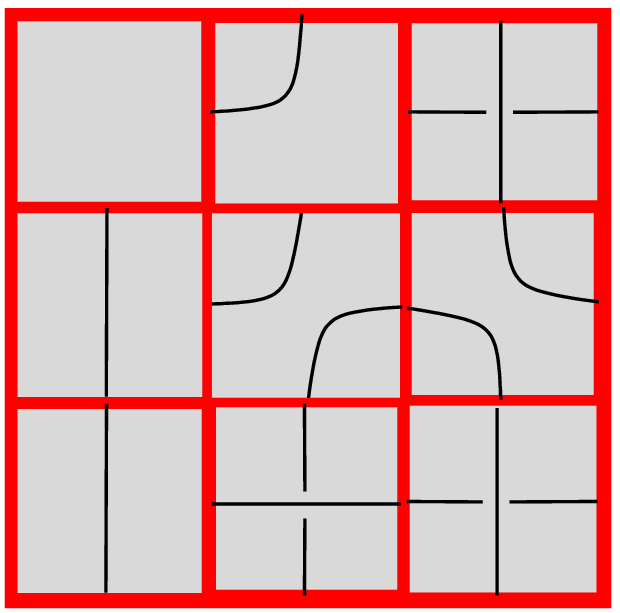}}
\caption{An example 3-mosaic }
\label{fig:examosaic}
\end{figure}

\begin{figure}[htb] \epsfysize = 0.7 in
\centerline{\epsffile{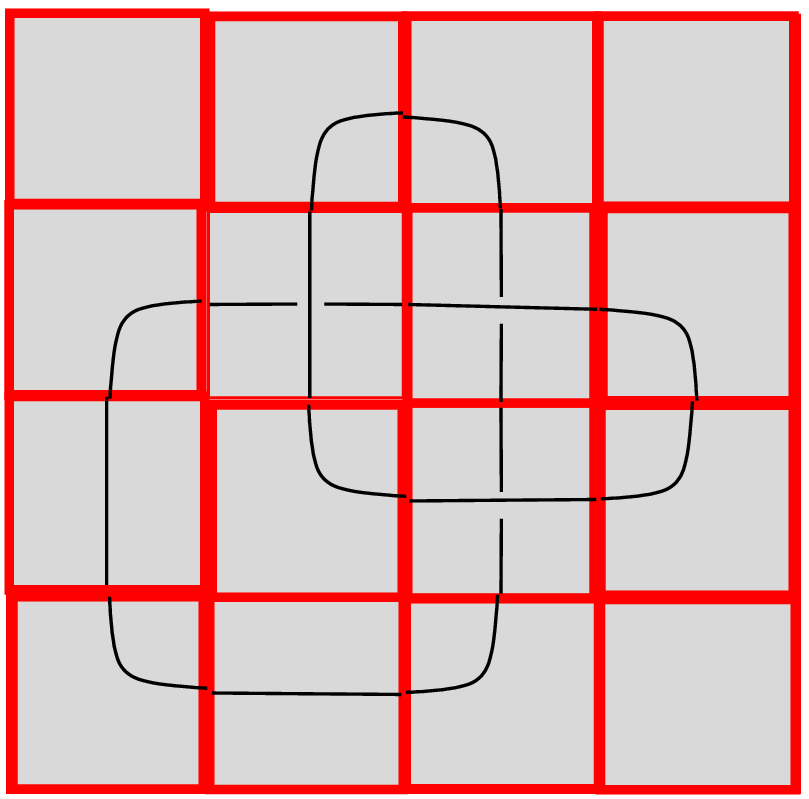}}
\caption{The trefoil as a mosaic knot }
\label{fig:trefoil}
\end{figure}

There are 11 mosaic tiles that are shown, up 90 degree rotations, in figure \ref{fig:tiles}. Each tile contains a tangle and has either $0$, $2$, or $4$ connection points. An \emph{ $n$-mosaic }is an $n \times n$ matrix of mosaic tiles.
 An example of a $3$-mosaic is shown in figure \ref{fig:examosaic}. 
Two adjacent tiles are said to be \emph{suitably connected} if their endpoints connect on the adjacent side. An \emph{$n$-mosaic link} is an $n \times n $ mosaic of suitably connected tiles. An $n$-mosaic knot in an $n$-mosaic link that contains exactly 1 component. An example of a $3$-mosaic knot is shown in figure 
\ref{fig:trefoil}. 

\begin{figure}[htb] \epsfysize = 6 in
\centerline{\epsffile{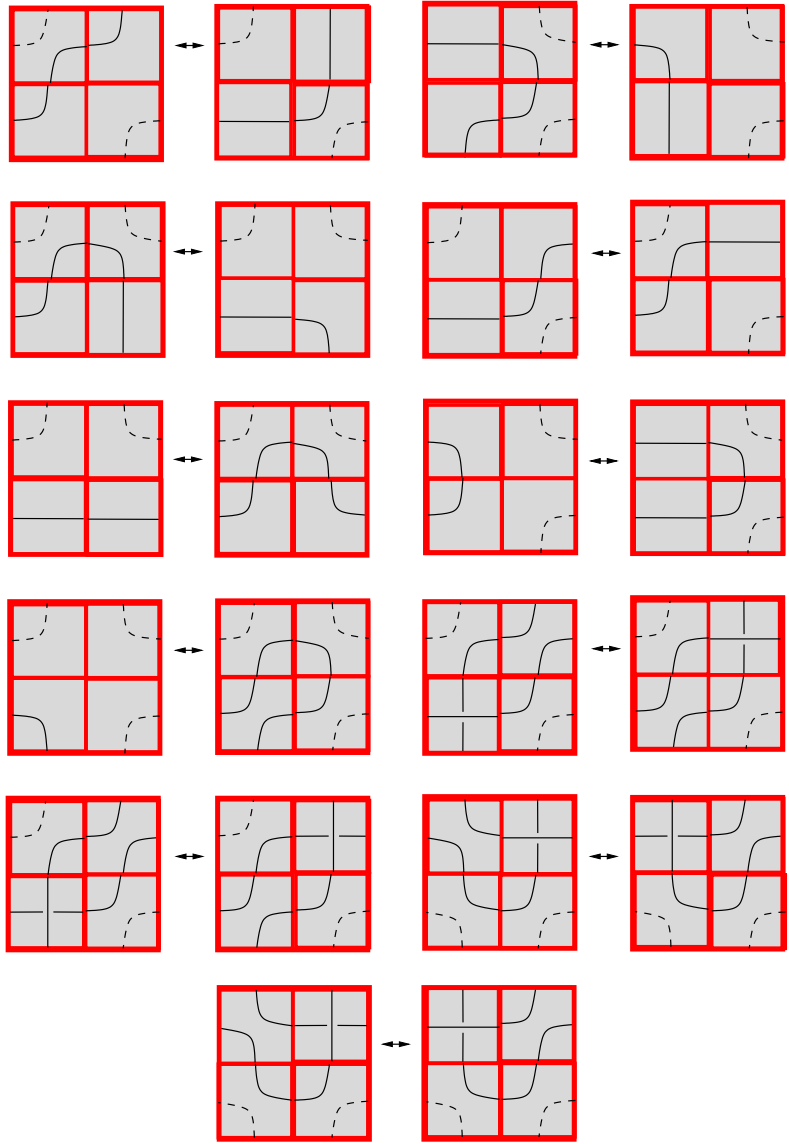}}
\caption{Mosaic planar istopy moves}
\label{fig:planar1}
\end{figure}

\begin{figure}[htb] \epsfysize = 2.5 in
\centerline{\epsffile{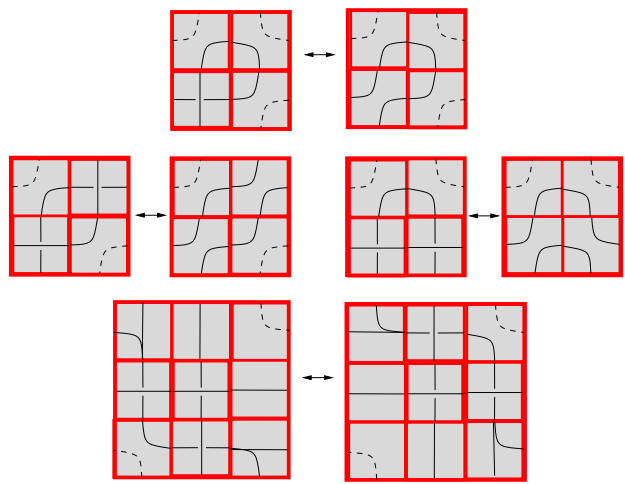}}
\caption{A subset of the mosaic Reidemeister moves}
\label{fig:rmoves}
\end{figure}

An $n$-mosaic link is equivalent to an $m$-mosaic link if they are related by a sequence of mosaic planar isotopy moves, mosaic Reidemeister moves, and increasing/decreasing the size of the mosaic grid.
The mosaic planar isotopy moves, up to ninety degree rotations, are shown in figure \ref{fig:planar1}. A subset of the mosaic Reidemeister moves are shown in figure \ref{fig:rmoves}. To obtain the full set of mosaic Reidemeister moves, we must all consider crossing changes and ninety degree rotations of the matrices of tiles. 
Equivalence classes of $n$-mosaic links are determined by sequences of mosaic planar isotopy moves, mosaic Reidemeister moves, and increasing or decreasing the size of the mosaic matrix.  A \emph{mosaic link} refers to an equivalence class of equivalent $n$-mosaic links.

\begin{rem}Takahito Kuriya proved that Mosaic knot equivalence classes
 are in 1-1 correspondence with Reidemeister move equivalence classes. We observe that this result does not extend to virtual mosaic knots. \end{rem}

\begin{figure}[htb] \epsfysize = 0.7 in
\centerline{\epsffile{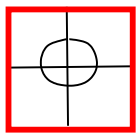}}
\caption{Virtual crossing tile }
\label{fig:vtile}
\end{figure}

\begin{rem}
We can extend this definition to virtual mosaic knots by introducing an additional tile, as shown in figure \ref{fig:vtile}. 
The definition of equivalence is expanded to include mosaic versions of the planar isotopy moves and the virtual Reidemeister moves. This extension will be discussed in more detail in section \ref{virtsection}.
\end{rem}

\begin{figure}[htb] \epsfysize = 0.5 in
\centerline{\epsffile{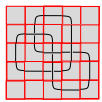}}
\caption{Mosaic figure Eight Knot}
\label{fig:fig8}
\end{figure}

The \emph{mosaic number} of a mosaic link is  an equivalence class is the smallest number $n$ for which a member of the equivalence class is an $n$-mosaic.  
The \emph{crossing number} of link is the smallest number of crossings in any diagram of the link, that is, the smallest number of crossings in any equivalent link diagram. For example, the
 crossing number of the trefoil is three and its mosaic number is four (see figure \ref{fig:trefoil}. The mosaic number of the figure eight knot is five and its crossing number is four as shown in figure \ref{fig:fig8}. This provokes several questions about the relationship between mosaic number and crossing numbers.

\section{An upper bound on crossing number} \label{up}

We first consider an upper bound on the crossing number.

\begin{thm} \label{above} Let $K$ be a knot or link with mosaic number $n$. Then the crossing number $c$ is bounded above:
\begin{equation}
c \leq (n-2)^2. 
\end{equation}
\end{thm}
\textbf{Proof:}
Observe that the edge and corner tiles of a mosaic cannot contain crossings if the mosaic is suitably connected. Since mosaic knots and links are suitably connected, only $(n-2)^2$ tiles can contain four connection points. Hence:
\begin{equation} 
c  \leq (n-2)^2. \qed
\end{equation}

\begin{figure}[htb] \epsfysize = 2 in
\centerline{\epsffile{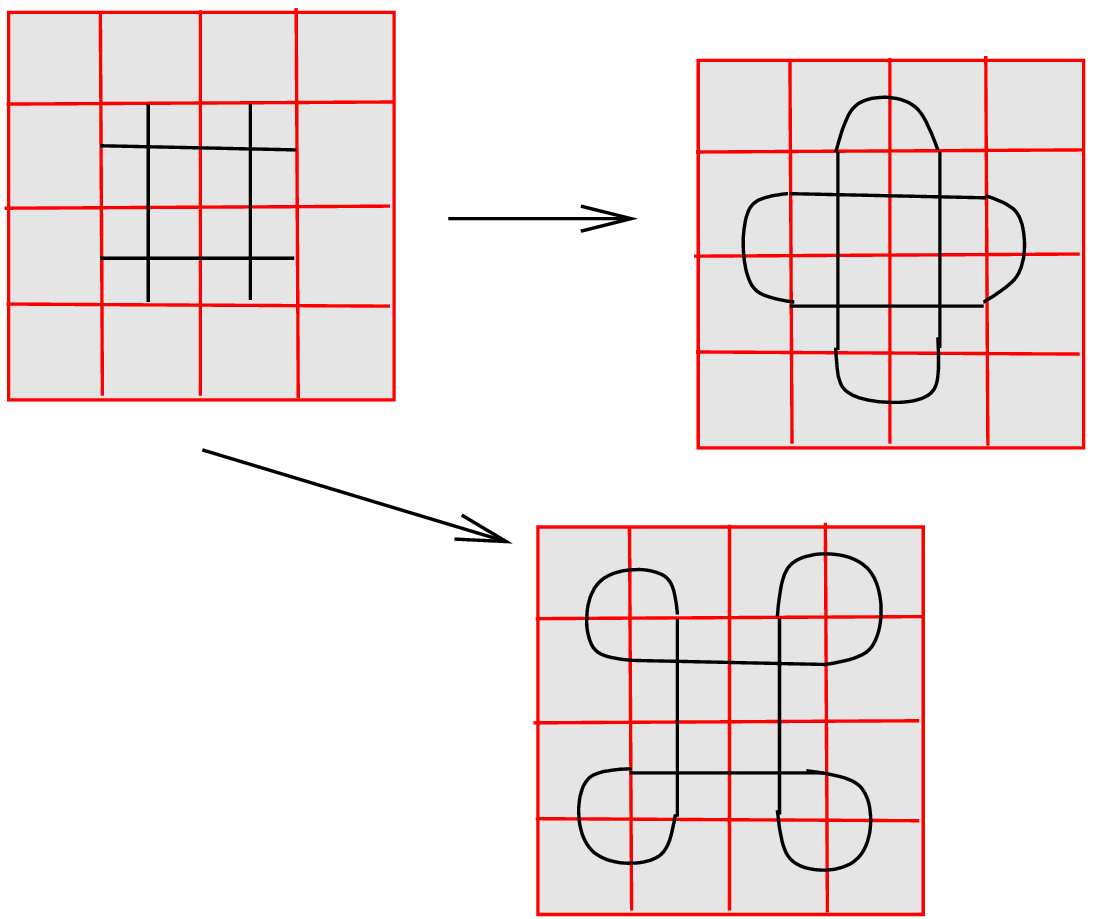}}
\caption{An even mosaic and possible closures}
\label{fig:even}
\end{figure}

\begin{cor}
Let $K$ be a knot with mosaic number $n$. If $n$ is even then $ c \leq (n-2)^2 -1 $ for $ n > 3$. 
\end{cor}
\textbf{Proof}
By Theorem \ref{above}, at most $ (n-2)^2 $ tiles contain crossings. If all 
$(n-2)^2 $ interior tiles contain crossings, then we either obtain a link or the unknot as shown in figure \ref{fig:even}. As a result, the mosaic can contain at most
$ (n-2)^2 -1 $ crossings and be a knot. \qed

\begin{figure}[htb] \epsfysize = 2 in
\centerline{\epsffile{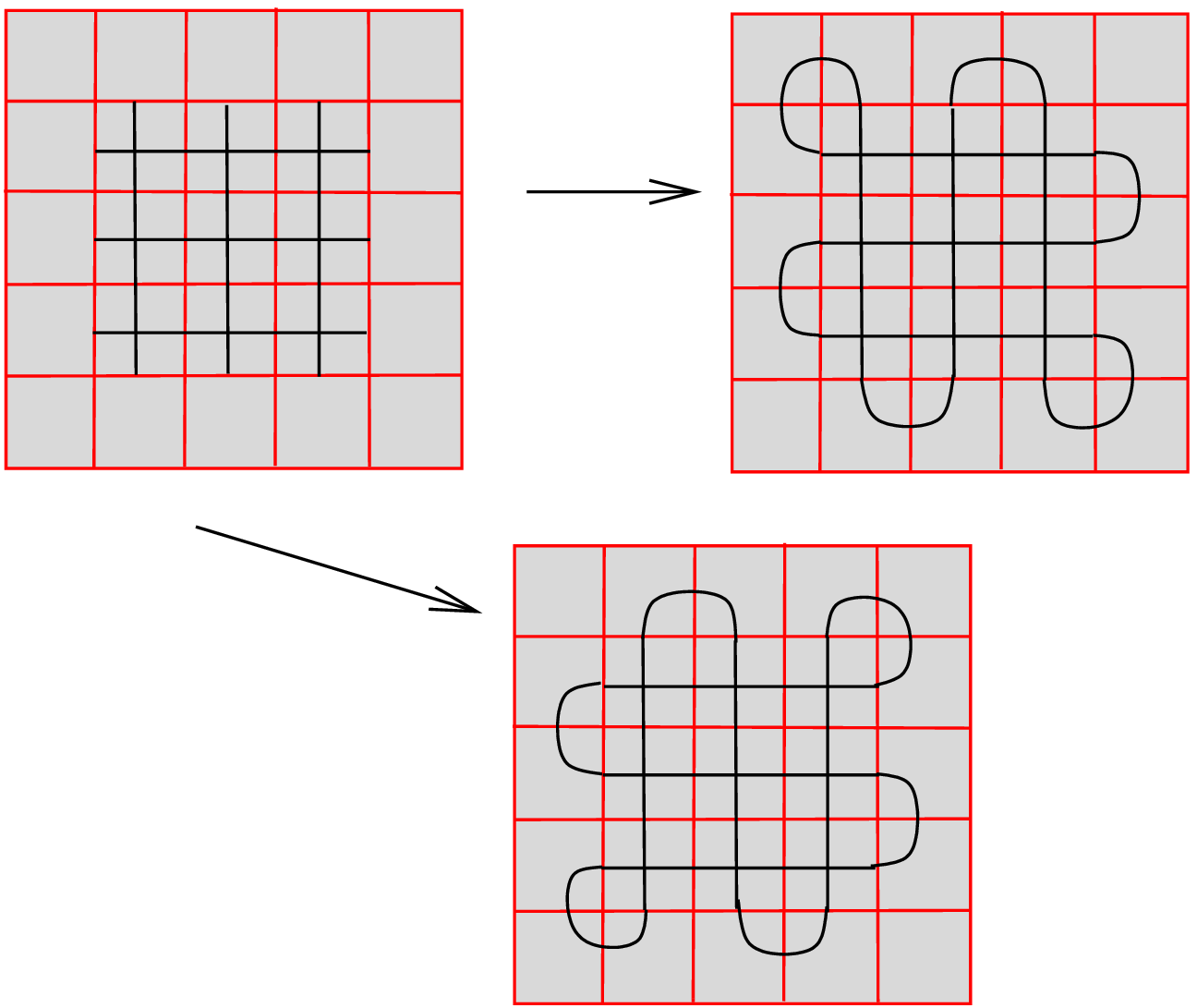}}
\caption{An odd mosaic and possible closures }
\label{fig:odd}
\end{figure}

This bound is sharp because the trefoil has mosaic number 4 and crossing number 3. See figure \ref{fig:trefoil}.

\begin{cor}
Let $K$ be a knot with mosaic number $n$. If $n$ is odd then $ c \leq (n-2)^2 -2 $ for $ n>3$.
\end{cor}
\textbf{Proof:} We first observe that by Theorem \ref{above}, the mosaic knot contains at most $ (n-2)^2$ crossings. If all $(n-2)^2$ tiles have crossings, then there are only two possible ways to obtain a suitably connected mosaic, as shown in figure \ref{fig:odd}. Both choices result in a diagram where it is possible to apply the Reidemeister I move to two corners. As a result, a knot with mosaic number $n$ has crossing number $c$ such that:
\begin{equation*}
c \leq (n-2)^2 -2. \qed 
\end{equation*}

\begin{figure}[htb] \epsfysize = 1 in
\centerline{\epsffile{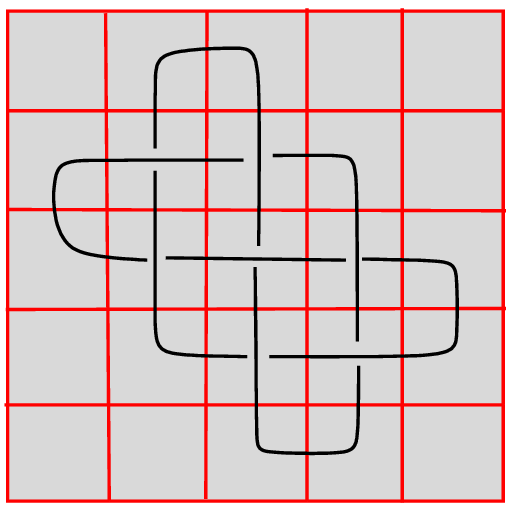}}
\caption{A knot with seven crossings and mosiac number 5 }
\label{fig:oddsharp}
\end{figure}

We again observe that this bound is sharp; consider the knot shown in \ref{fig:oddsharp}. This knot has crossing number seven and mosaic number five.

\section{A rough lower bound on mosaic number} \label{lower}

\begin{thm}Let $K$ be link with crossing number $c$ then the mosaic number $n$ is bounded below by:
\begin{equation}
\lceil \sqrt{c}  \rceil +2 \leq n.
\end{equation}
\end{thm}
\textbf{Proof:}
Let $K$ be a mosaic link with crossing number $c$ and let $n$ be the mosaic number. Now, $K$ can contain at most $(n-2)^2 $ crossings. Hence $ c \leq (n-2)^2 $. Rewriting, we obtain:
\begin{equation*}
\lceil \sqrt{c} \rceil +2 \leq n .\qed
\end{equation*}

\section{An upper bound on  the mosaic number}

\begin{figure}[htb] \epsfysize = 4 in
\centerline{\epsffile{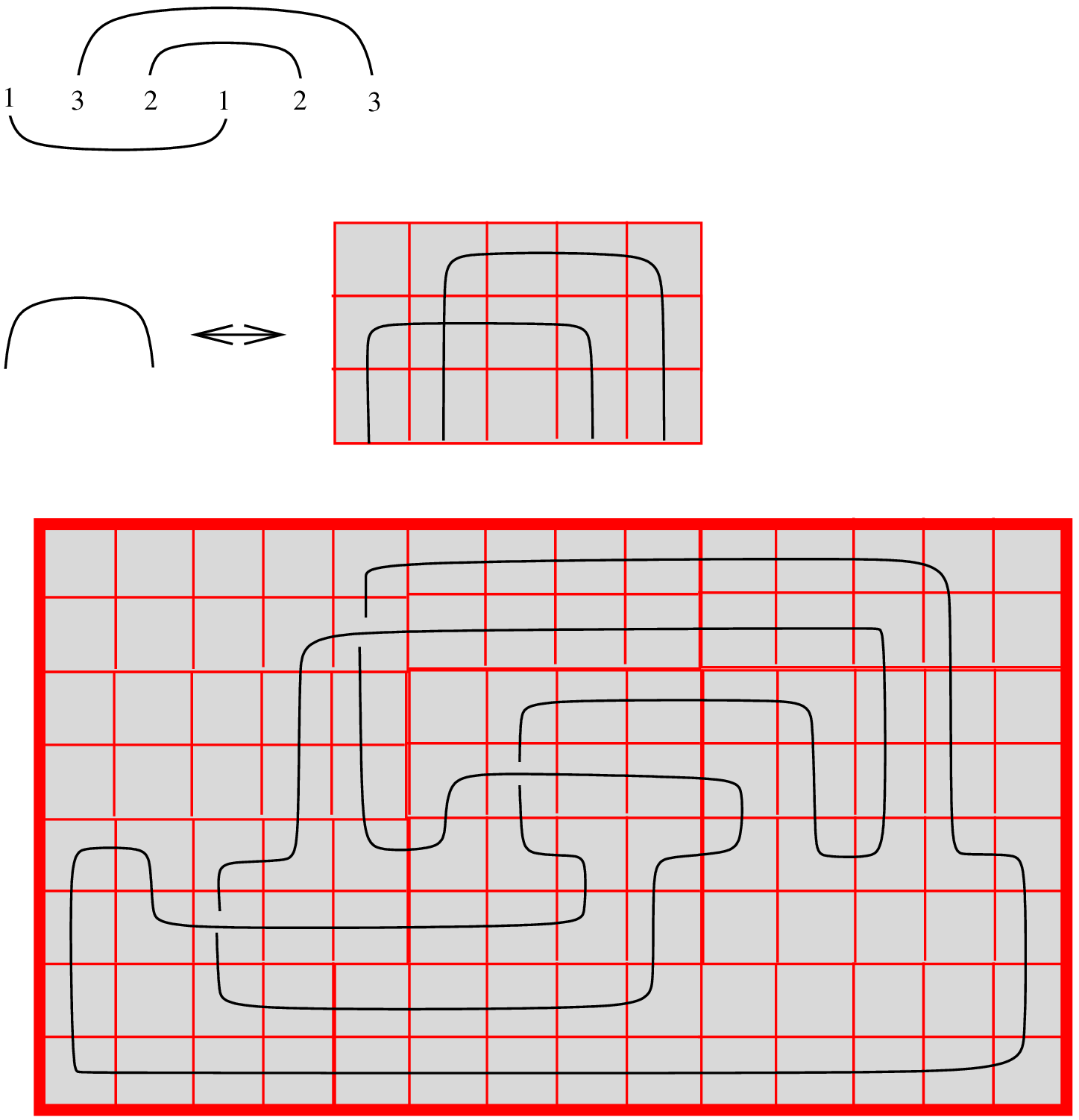}}
\caption{Constructing a mosaic knot }
\label{fig:process}
\end{figure}

\begin{thm} \label{upper} Let $K$ be a knot diagram with crossing number $c$. Then the mosaic number $n$ can be bounded above as follows:
\begin{equation}
4c+2 \geq n.
\end{equation}
\end{thm}
\textbf{Proof:}
From the link diagram $K$ with crossing number $c$, we can construct the Gauss code by fixing a point on each component of the link and a direction of orientation. Number each crossing in the link  and traverse each component from the fixed point in the direction of orientation. Record the number  of each crossing as it is traversed, in order to obtain a list of labels of length $2c$ where each label is recorded twice. Following the algorith given in \cite{kvirt}, we then invert the list of numbers, following the order of the numbers. That is, invert the list between the two labels of \textit{1} and then invert the list between the labels of \textit{2} and so on. From this, we can reconstruct the link by inserting 2-2 tangles with a single crossing as shown in figure \ref{fig:process}. For each crossing, we need 4 tiles and for the closure of the knot diagram, we require 2 additional tiles across the width of the mosaic. Vertically, each crossing requires a two rows of tiles, plus a row tiles for the center and a row of tiles for the closure. As a result:
\begin{equation*}
4c+2 \geq n. \qed
\end{equation*}

\begin{cor}
Let $K$ be a knot with crossing number $c$ then the mosaic number $n$ is bounded above: $ 4c+2 \geq n$.
\end{cor}
\textbf{Proof:} Let $K$ be a knot with crossing number $c$. Then $K$ has a diagram with $c$ crossings. \qed

Based on experimental evidence, we conjecture that:
\begin{conj}The mosaic number $n$ is bounded above by $c+2$. 
That is, $ n \leq c+2 $. 
\end{conj}

\section{An extension of these bounds to virtual mosaics} \label{virtsection}

\begin{figure}[htb] \epsfysize = 2 in
\centerline{\epsffile{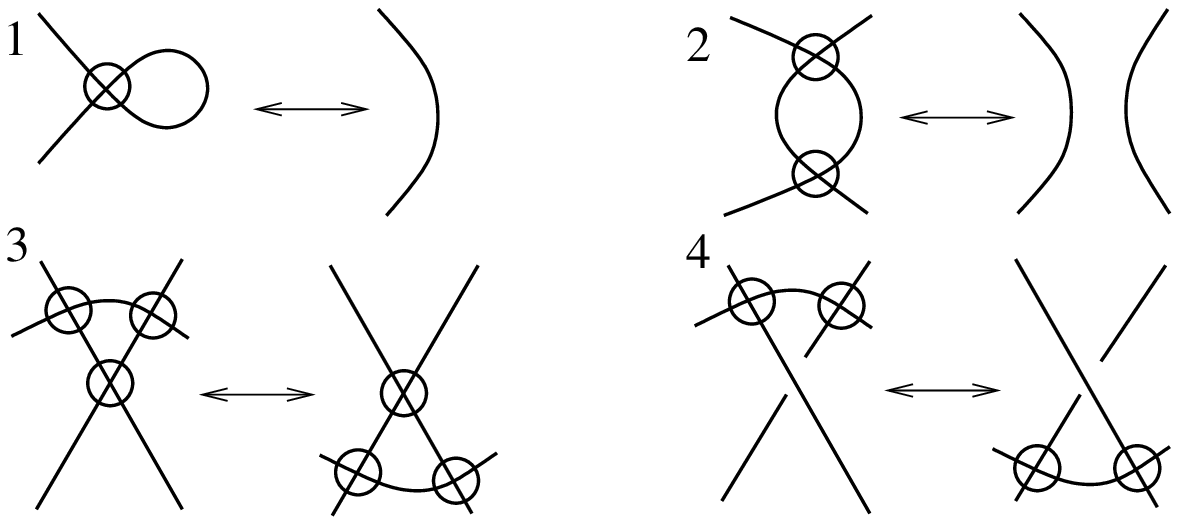}}
\caption{The virtual Reidemeister moves }
\label{fig:vrmoves}
\end{figure}

Virtul link diagrams are decorated immersions of $n$ copies of $S^1 $ into the plane with two types of crossings: classical (indicated by over/under markings) and virtual (indicated by a solid encircled X) \cite{kvirt}. 
We recall that two virtual link diagrams are said to be \emph{virtually equivalent} if one can be transformed into the other by a sequence of virtual and classical Reidemeister moves. A \emph{virtual link} is an equivalence class of virtual link diagrams.
The virtual Reidemeister moves are shown in figure \ref{fig:vrmoves}.

\begin{figure}[htb] \epsfysize = 3 in
\centerline{\epsffile{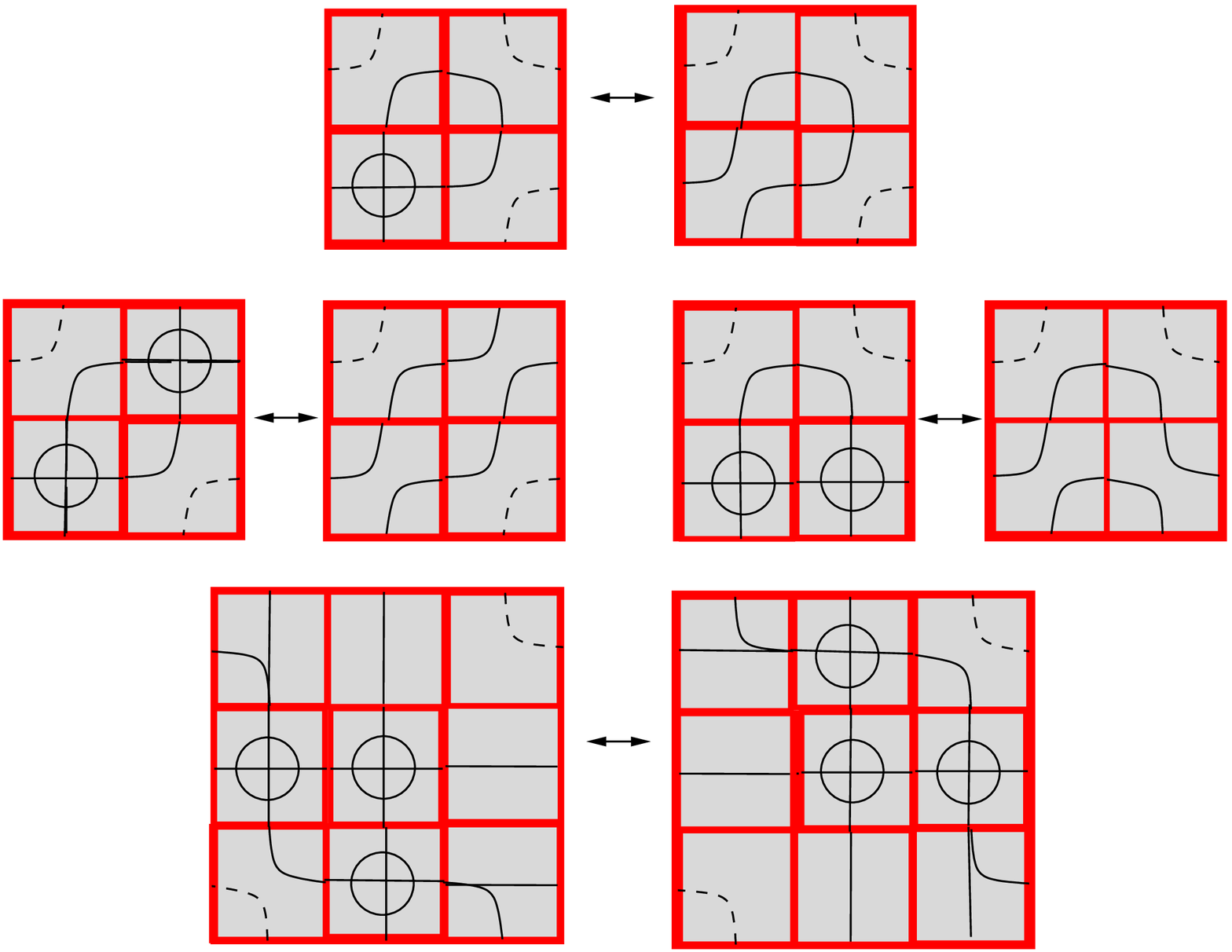}}
\caption{Mosaic virtual Reidemeister moves }
\label{fig:mvrmoves}
\end{figure}

We can extend mosaic knot theory to virtual knots by incorporating an additional type of tile: a virtual crossing tile and extending the mosaic Reidemeister moves to include the virtual Reidemeister moves. We define a \emph{virtual $n$-mosaic} to be a $n \times n$ mosaic of tiles that possibly includes the virtual crossing tile. A \emph{virtual $n$-mosaic link } is a suitably connected virtual $n$-mosaic.  Mosaic versions of some of the virtual Reidemesiter moves are shown in figure \ref{fig:mvrmoves}.

A virtual $n$-mosaic link diagram  is equivalent to a virtual $m$-mosaic didagram is one may be transformed in the the other by a sequence of mosaic planar isotopy moves, mosaic Reidemeister moves and virtual mosaic Reidemeister moves, as well as increasing or decreasing the size of the mosaic. A \emph{virtual mosaic link} is an equivalence classe of virtual $n$ mosaic link diagrams.
\begin{rem}
For a complete description of virtual $n$-mosaics, see \cite{irina}. 
\end{rem} 

Let $k$ denote the \emph{virtual crossing number} of a virtual link. The virtual crossing number is the minimum number of virtual crossings in any diagram of the virtual link. We can extend the result of Theorems \ref{up} and \ref{lower} to the virtual case. 

\begin{thm}Let $K$ be a virtual mosaic link with mosaic number $n$, crossing number $c$, and virtual crossing number $k$. Then:
\begin{equation}
c+k \leq (n-2)^2.
\end{equation}
\end{thm}
\textbf{Proof:} Analogous to the proof of Theorem \ref{upper}.

\section{Other questions on mosaic knots}
There are a variety of interesting questions about mosaic knots that can be approached from a combinatorial viewpoint. These questions will be considered in future work.
\begin{itemize}
\item What is the maximum number of components that a mosaic link with mosaic number $n$ can contain?
\item Can Kuriya's result be extended to the virtual case? 
\item What is the maximum size diagram needed to convert an $n$-mosaic into a $m$-mosaic? Can this number be bounded?
\end{itemize} 

\textbf{Acknowledgement:} This paper details the work completed during a Research Experience for Undergraduates at McKendree University during summer 2009. This program  was a Mathematical Association of America activity funded by National Security Agency (grant H98230-06-1-0156) and the National Science Foundation (grant DMS-0552763).

\end{document}